\documentclass[11pt]{article}

\usepackage{amsmath}
\usepackage{amssymb}
\usepackage{indentfirst}
\usepackage{graphics} 
\usepackage{color}

\makeatletter

\@addtoreset{equation}{section}
\makeatother
\def\<{\langle}
\def\>{\rangle}

\newtheorem{theo}{Theorem}[section]
\newtheorem{rem}{Remark}[section]

\newtheorem{cor}{Corollary}[section]
\makeatletter
   
   \@addtoreset{equation}{section}
\makeatother

\begin{document}
\title{\bf Some attempts on $L^2$ boundedness for $1$-D \\wave equations with time variable coefficients}
\author{Ryo Ikehata\thanks{ikehatar@hiroshima-u.ac.jp} \\ {\small Department of Mathematics, Division of Educational Sciences}\\ {\small Graduate School of Humanities and Social Sciences} \\ {\small Hiroshima University} \\ {\small Higashi-Hiroshima 739-8524, Japan} \\}
\date{}
\maketitle
\begin{abstract}
We consider the $L^{2}$-boundedness of the solution itself of the Cauchy problem for wave equations with time-dependent wave speeds. We treat it in the one-dimensional Euclidean space ${\bf R}$. To study these, we adopt a simple multiplier method by using a special property equiped with the one dimensional space. 
\end{abstract}
\section{Introduction}
\footnote[0]{Keywords and Phrases: wave equation; $1$-D space; time-variable coefficient; $L^{2}$-boundedness; multiplier method.}
\footnote[0]{2020 Mathematics Subject Classification. Primary 35L05; Secondary 35L10, 35B45, 35B40.}

We consider the Cauchy problem
 for wave equations with time dependent wavespeeds in one-dimensional Euclidean space ${\bf R}$
\begin{equation}
u_{tt}(t,x) - a(t)^{2}u_{xx}(t,x) = 0,\ \ \ (t,x)\in (0,\infty)\times {\bf R},\label{eqn}
\end{equation}
\begin{equation}
u(0,x)= u_{0}(x),\ \ u_{t}(0,x)= u_{1}(x),\ \ \ x\in {\bf R},\label{initial}
\end{equation}
where the initial data $[u_{0},u_{1}]$ are taken from the test function space (for simplicity)
\[u_{0} \in C_{0}^{\infty}({\bf R}),\quad u_{1} \in C_{0}^{\infty}({\bf R}),\]
$a \in {\rm C}^{1}([0,\infty))$ and we denote
\[u_{t}=\frac{\partial u}{\partial t},\quad u_{tt}=\frac{\partial^2 u}{\partial t^2},\quad u_{xx}=\frac{\partial^2 u}{\partial x^2}.\]

Throughout this paper, $\| \cdot\|$ stands for the usual $L^2({\bf R})$-norm. 
The total energy $E_{u}(t)$ of the solution $u(t,x)$ to problem (1.1) is defined by
\begin{equation}
E_{u}(t)=\frac{1}{2}(\| u_t(t,\cdot)\|^2+a(t)^{2}\| u_x(t,\cdot)\|^2).
\end{equation}

We shall impose the following two assumptions on $a(t)$:\\
({\bf A.1})\,$a \in {\rm C}^{1}([0,\infty))$ and $a(t) > 0$ for all $t \geq 0$, and $a_{m} := \sup\{a(t)\,:\,t \geq 0\} < +\infty$,\\
({\bf A.2})\,$a'(t) \geq 0$ for all $t \geq 0$.\\
\noindent
Under these conditions, it is known that the problem (1.1)-(1.2) admits a unique strong solution $u \in {\rm C}([0,\infty);H^{2}({\bf R})) \cap {\rm C}^{1}([0,\infty);H^{1}({\bf R}))\cap {\rm C}^{2}([0,\infty);L^{2}({\bf R}))$, which has finite speed propagation of waves with its propagation speed $a_{m} > 0$ (cf. \cite{ikawa}). This is a sufficient class to deploy the multiplier method. In fact, the solution is much smoother than is convenient.\\

The purpose of the present paper is to consider whether the $L^{2}$-boundedness of the solution itself for problem (1.1)-(1.2) can be observed or not in the one-dimensional case. The problem itself is never trivial in the sense that one has no Hardy or Poincar\'e inequality. Furthermore, one can imagine that time-dependent coefficients $a(t)$ can be an obstacle in many ways in partial integration.

There have been many interesting papers published on the estimate of the total energy and asymptotic behavior of solutions of the wave equation, not only in the case of constant coefficients, but also in the case of time and space variable dependence (see \cite{AR, C, E, H, HW, Ma, RS, Y} and the references therein). Furthermore, a series of sharp research results have been published on the $L^p$-$L^q$-estimate of the solution itself, starting with Strichartz \cite{S} (see e.g. \cite{BS, B, M, P, RY, W}). However, when we look at the $L^2$ estimate of the solution itself, there seems to be a little gap from the viewpoint of bestness. In such a situation, the author \cite{JHDE-ike} published the following results for the case $a(t) = 1$ (constant coefficient) in problems (1.1)-(1.2):
\begin{equation}\label{1.4}
\int_{-\infty}^{\infty}u_{1}(x)dx \ne 0 \Rightarrow \Vert u(t,\cdot)\Vert \sim \sqrt{t},\quad (t \to \infty).
\end{equation}
This shows a growth estimate of the $L^{2}$-norm of the solution itself to problem (1.1)-(1.2) with $a(t) = 1$. Thus, when one wishes to observe the $L^2$ boundedness of the solution itself, it will be necessary to consider the general time-variable coefficient problem (1.1) by factoring in the information that the zero-order moment of the initial velocity may or may not vanish. Herein lies the difficulty of deriving the solution's own $L^2$ estimate in the case of general variable coefficients. Furthermore, with variable coefficients, it is not possible to describe and evaluate the solution "explicitly" using the Fourier transform as in the method of \cite{JHDE-ike}, and this makes us imagine a very hopeless prediction for deriving the $L^2$ estimate formula from "under" the solution itself. Therefore, in order to construct a "general theory" involving constant coefficients to watch the $L^2$ boundedness of the solution itself, the condition that the zero-order moment of the initial velocity vanishes, can never be removed.\\
It's time to introduce the results.
%
To state our results, we define
\[v_{1}(x) := \int_{-\infty}^{x}u_{1}(y)dy.\]
Our result then can be stated as follows. 
\begin{theo} \label{th1}
Suppose {\rm ({\bf A.1})}, {\rm ({\bf A.2})} and assume $v_{1} \in L^{2}({\bf R})$. Then, the corresponding solution $u(t,x)$ to problem {\rm (1.1)-(1.2)} with initial data $[u_{0},u_{1}] \in C_{0}^{\infty}({\bf R})\times C_{0}^{\infty}({\bf R})$ satisfies
\[\Vert u(t,\cdot)\Vert^{2} \leq I_{0}^{2}a(0)^{-2},\quad (t \geq 0)\]
with a constant $I_0 \geq 0$ defined by 
\[I_{0} := \left( \Vert v_{1}\Vert^{2} + a(0)^{2}\Vert u_{0}\Vert^{2} \right)^{1/2}.\]
\end{theo}
\begin{rem}{\rm In this Theorem \ref{th1}, it is essentially imposing a new condition on $u_1$ through $v_{1} \in L^{2}({\bf R})$. For example, if $u_{1} \in C_{0}^{\infty}({\bf R})$ is odd about the origin, then we see that $\int_{{\bf R}}u_{1}(x)dx = 0$, and in this case $v_{1} \in C_{0}^{\infty}({\bf R})$. Thus one has $v_{1} \in L^{2}({\bf R})$. The condition $v_{1} \in L^{2}({\bf R})$ on the initial velocity $u_{1}$ of the theorem makes sense.}
\end{rem}
\noindent
{\bf Example 1.}\, We can find an appropriate function $a(t)$ satisfying {\rm ({\bf A.1})} and {\rm ({\bf A.2})}. \\
\[a(t) = \left\{
  \begin{array}{ll}
   \displaystyle  {1+e^{-1/t}},& \qquad t > 0,\\[0.2cm]
   \\
   \displaystyle {1},&
   \qquad t = 0,
   \end{array} \right. \]
\noindent
From the proof of Theorem \ref{th1} we see that the condition {\bf (A.2)} can be replaced by the following one:\\
({\bf A.3})\,$a'(t) \leq 0$ for all $t \geq 0$, and $A_{0} := \inf\{a(t)\,:\,t \geq 0\} > 0$.\\
Then, one can also derive the following corollary. A monotone "decreasing" function $a(t)$ can be also covered in our theory.
\begin{cor} \label{th3}
Suppose {\rm ({\bf A.1})} and {\rm ({\bf A.3})} and assume $v_{1} \in L^{2}({\bf R})$. Then, the corresponding solution $u(t,x)$ to problem {\rm (1.1)-(1.2)} with initial data $[u_{0},u_{1}] \in C_{0}^{\infty}({\bf R})\times C_{0}^{\infty}({\bf R}))$ satisfies
\[\Vert u(t,\cdot)\Vert^{2} \leq I_{0}^{2},\quad (t \geq 0).\]
\end{cor}
\noindent
{\bf Example 2.}\, We can also choose $a(t) := 1+e^{-t}$, and/or $a(t) := \frac{2+t}{1+t}$. Then, the statement of Corollary \ref{th3} implies
\[\Vert u(t,\cdot)\Vert^{2} \leq I_{0}^{2},\quad (t \geq 0).\]
\begin{rem}{\rm The conditions $\sup\{a(t)\,:\,t \geq 0\} < +\infty$ and $\inf\{a(t)\,:\,t \geq 0\} > 0$ assumed in {\bf (A.2)} and {\bf (A.3)} express a propagation speed of the wave, and the ellipticity of the operator $u \mapsto a(t)^{2}\partial_{xx}u$, respectively.}
\end{rem}
By modifying the proof of Theorem \ref{th1} one can also present another version of the result. For this purpose we set one more assumption. \\
({\bf A.4})\,$A_{0} > 0$, and $a' \in L^{1}(0,\infty)$,\\
where $A_{0}$ is the constant already defined in {\bf (A.3)}. Then, one can derive one more corollary. 
\begin{cor} \label{th4}
Suppose {\rm ({\bf A.1})} and {\rm ({\bf A.4})} and assume $v_{1} \in L^{2}({\bf R})$. Then, the corresponding solution $u(t,x)$ to problem {\rm (1.1)-(1.2)} with initial data $[u_{0},u_{1}] \in C_{0}^{\infty}({\bf R})\times C_{0}^{\infty}({\bf R}))$ satisfies
\[\Vert u(t,\cdot)\Vert^{2} \leq \frac{1}{A_{0}^{2}}I_{0}^{2}e^{\frac{2}{A_{0}}\int_{0}^{\infty}\vert a'(s)\vert ds},\quad (t \geq 0).\]
\end{cor}
\noindent
{\bf Example 3.}\, We can present an additional oscillating example:
\[a(t) := 2 + \frac{\sin t}{(1+t)^{2}}.\]

The above theorem and examples include the case of constant coefficients $a(t) = 1$. Therefore, what is still a concern is the fear that some moment conditions for the initial velocity $u_{1}$ may contradict (1.4). Let us discuss this situation below.
Suppose $u_{1} \in C_{0}^{\infty}({\bf R})$. Then, there is a large number $L > 0$ such that one can assume that ${\rm supp}\,u_{1} \subset [-L,L]$. Since
\[v_{1}(x) = \int_{-\infty}^{x}u_{1}(y)dy,\]
if $x > 2L$, then one sees that
\[v_{1}(x) = \int_{-L}^{L}u_{1}(y)dy = \int_{-\infty}^{\infty}u_{1}(y)dy =: c_{0}\quad (\forall x > 2L).\]
Assume for the moment that $c_{0} \ne 0$. Then, it follws that
\[\Vert v_{1}\Vert^{2} = \int_{-\infty}^{2L}\vert v_{1}(x)\vert^{2}dx + \int_{2L}^{\infty}\vert v_{1}(x)\vert^{2}dx \geq \int_{2L}^{\infty}\vert v_{1}(x)\vert^{2}dx = \int_{2L}^{\infty}c_{0}^{2}dx = \infty,\]
which shows a contradiction to the assumption that $v_{1} \in L^{2}({\bf R})$ in Theorem 1.1 and Corollary 1.1. Thus it must be $c_{0} = 0$. Note that $u_{1} \in C_{0}^{\infty}({\bf R})$ and $v_{1} \in L^{2}({\bf R})$ impliy $v_{1} \in C_{0}^{\infty}({\bf R})$ and $v_{1}(\infty) = 0$. Thus, under the conditions of Theorem \ref{th1}, there is no contradiction with \eqref{1.4} because of the
additional assumption that the zero-order moment of the initial velocity vanishes. Conversely, if the zero-order moment of the initial velocity does not vanish, then the augmentation property as in \eqref{1.4} may still be derived, but this remains
unresolved at this time. Note that the theorem took the initial value class to be the test function, but this is not the essence of the theorem. This is off the topic of this issue, though.
\begin{rem}{\rm With our method, it is difficult to handle high dimensional cases now and in our time. Also, by analogy with the results in \cite{JHDE-ike} for the constant coefficients, the moment vanishing condition for initial velocity $u_{1}$ would be required even in two dimensions, but that condition would not be necessary for three or more dimensions. The $L^2$-bounddedness is more likely to be used for three or more dimensions.}
\end{rem}

\begin{rem}{\rm From the above discussion, we expect the following in the $1$-D case (and probably in the $2$-D case as well): under the condition such that
\[\int_{{\bf R}}u_{1}(x)dx \ne 0,\]
then
\[\lim_{t \to \infty}\Vert u(t,\cdot)\Vert = +\infty\]
with some growth rate. This is still open in the time variable coefficient case (cf. \eqref{1.4}).}
\end{rem}

The remainder of this paper is organized as follows. In Section 2, we shall prove Theorem 1.1 and Corollaries 1.1 and 1.2.

\section{Proof of results.}

In this section, let us prove Theorem \ref{th1} and then Corollary \ref{th3} by using multiplier method inspired from the idea in \cite{Z}. Incidentally, it seems that the method used in \cite{Z} itself is partially inspired from an idea developed in \cite{IM} as can be observed from their proofs.\\

{\it Proof of Theorem \ref{th1}.}\\
First of all, define a function $v(t,x)$ by
\[v(t,x) := \int_{-\infty}^{x}u(t,y)dy.\]
Note that the function $v(t,x)$ is well-defined because of finite speed propagation of waves. Furthermore, for each $j = 0,1$ we set
\[v_{j}(x) := \int_{-\infty}^{x}u_{j}(y)dy.\]
Incidentally, we see that $v_{j} \in C^{\infty}({\bf R})$ ($j = 0,1$). Then, the function $v(t,x)$ satisfies
\begin{equation}
v_{tt}(t,x) - a(t)^{2}v_{xx}(t,x) = 0,\ \ \ (t,x)\in (0,\infty)\times {\bf R},\label{eqn-v}
\end{equation}
\begin{equation}
v(0,x) = v_{0}(x),\ \ v_{t}(0,x) = v_{1}(x),\ \ \ x\in {\bf R},\label{initial-v}
\end{equation}
where one has just used the fact that
\[\lim_{y \to -\infty}u_{y}(t,y) = 0\quad (t \geq 0).\]
Multiplying both sides of \eqref{eqn-v} by $v_{t}$, and integration by parts yield the equality:
\[\frac{d}{dt}E_{v}(t) = a(t)a'(t)\Vert v_{x}(t,\cdot)\Vert^{2},\]
where one has used the fact that $v_{x}(t,x) = u(t,x) = 0$ for large $\vert x\vert \gg 1$ and each $t \geq 0$.
Integrating the above equality on $[0,t]$ it follows that
\begin{equation}\label{ike-1}
E_{v}(t) = E_{v}(0) + \int_{0}^{t}a(s)a'(s)\Vert v_{x}(s,\cdot)\Vert^{2}ds
\end{equation}
\begin{equation}\label{ike-2}
= E_{v}(0) + 2\int_{0}^{t}\frac{a'(s)}{a(s)}\left(\frac{1}{2}a(s)^{2}\Vert v_{x}(s,\cdot)\Vert^{2}\right)ds
\end{equation}
\[\leq E_{v}(0) + 2\int_{0}^{t}\frac{a'(s)}{a(s)}E_{v}(s)ds,\]
where we have just used the assumption {\bf (A.2)}. Thus, by using the Gronwall inequality one has
\begin{equation}\label{i1}
E_{v}(t) \leq E_{v}(0)e^{2\int_{0}^{t}\frac{a'(s)}{a(s)}ds}.
\end{equation}
Therefore, from the definition of the total energy and \eqref{i1} one can get the inequality:
\[
\Vert v_{x}(t,\cdot)\Vert^{2} \leq 2E_{v}(0)a(t)^{-2}e^{2W(t)},
\]
where 
$$W(t) := \int_{0}^{t}\frac{a'(s)}{a(s)}ds = \int_{0}^{t}\frac{d}{ds}\log a(s) ds = \log \frac{a(t)}{a(0)}.$$
Thus one has
\[
\Vert v_{x}(t,\cdot)\Vert^{2} \leq 2E_{v}(0)a(t)^{-2}\left( \frac{a(t)}{a(0)}\right)^{2} = 2E_{v}(0)a(0)^{-2}.
\]
Now, since $v_{x}(t,x) = u(t,x)$ (This is a crucial idea in \cite{IM} and \cite{Z}), we can arrive at the desired estimate:
\[
\Vert u(t,\cdot)\Vert^{2} \leq 2E_{v}(0)a(0)^{-2} \quad (t \geq 0).
\]  
Set
\begin{equation}\label{i2}
I_{0} := \left( 2E_{v}(0) \right)^{1/2} = \left(\Vert v_{1}\Vert^{2} + a(0)^{2}\Vert u_{0}\Vert^{2} \right)^{1/2}.
\end{equation}
These imply the desired statement of Theorem \ref{th1}. 
\hfill
$\Box$
\\
\begin{rem}{\rm 
The function $v(t,x) := \int_{-\infty}^{x}u(t,y)dy$ used in the proof above can be considered as
\[v(t,x) \sim v(t,\infty) \quad (x \to \infty),\]
and
\[v(t,\infty) = \int_{{\bf R}}u(t,y)dy = \left. \int_{{\bf R}}e^{-iy\xi}u(t,y)dy\right\vert_{\xi = 0} = {\cal F}(u(t,\cdot))(0) = \hat{u}(t,0),\]
where ${\cal F}$ denotes the Fourier transform. Thus, the function $v(t,x)$ is equal to $\hat{u}(t,0)$ as $x \to \infty$ approximately. The estimate for $v(t,x)$ may correspond to the low frequency estimate near $\xi = 0$ of $\hat{u}(t,\xi)$. }
\end{rem}


\vspace{0.2cm}
{\it Proof of Corollary \ref{th3}.}\\
It follows from \eqref{ike-1} of the proof of Theorem \ref{th1} and the assumption {\bf (A.3)} we see that
\[E_{v}(t) \leq E_{v}(0).\]
The other part of proof is similar.
\hfill
$\Box$
\\

\vspace{0.3cm}
{\it Proof of Corollary \ref{th4}.}\\
It follows from \eqref{ike-2} of the proof of Theorem \ref{th1} and the assumption {\bf (A.4)} we see that
\[E_{v}(t) \leq E_{v}(0) + 2\int_{0}^{t}\frac{\vert a'(s) \vert}{a(s)}E_{v}(s)ds \leq E_{v}(0) + \frac{2}{A_{0}}\int_{0}^{t}\vert a'(s) \vert E_{v}(s)ds.\]
The use of the Gronwall inequality yields the desired estimate similarly.
\hfill
$\Box$


\vspace{0.5cm}
\noindent{\em Acknowledgement.}
\smallskip
The work of the author (R. IKEHATA) was supported in part by Grant-in-Aid for Scientific Research (C) 22540193 of JSPS.


\end{document}